\documentclass{amsart}
\usepackage{amsmath,amscd,amsthm,amsfonts,amssymb}

\theoremstyle{plain}
\newtheorem{theorem}                 {Theorem}      [section]
\newtheorem{proposition}  [theorem]  {Proposition}

\theoremstyle{definition}
\newtheorem{example}      [theorem]  {Example}
\newtheorem{remark}       [theorem]  {Remark}
\newtheorem{definition}   [theorem]  {Definition}

\numberwithin{equation}{section}

\def \theo-intro#1#2 {\vskip .25cm\noindent{\bf Theorem #1\ }{\it #2}}

\def \zn{\mathbb Z}

\def \rn{\mathbb R}

\def \F{\mathcal F}
\def \H{\mathcal H}

\def \V{\mathcal V}

\def \ip #1#2{\langle #1,#2 \rangle}

\def \lb#1#2{[#1,#2]}

\def \g{\mathfrak{g}}

\def \n{\mathfrak{n}}

\def \s{\mathfrak{s}}

\DeclareMathOperator{\ad}{ad}

\def \so#1{\mathfrak{so}(#1)}

\def \co#1{\mathfrak{co}(#1)}

\def \nab#1#2{\hbox{$\nabla$\kern -.3em\lower 1.0 ex
    \hbox{$#1$}\kern -.1 em {$#2$}}}

\begin{document}
\baselineskip 22pt \larger

\allowdisplaybreaks

\title{Harmonic morphisms from \\ four-dimensional Lie groups}

\author{Sigmundur Gudmundsson}
\author{Martin Svensson}



\address
{Department of Mathematics, Faculty of Science, Lund University,
Box 118, S-221 00 Lund, Sweden}
\email{Sigmundur.Gudmundsson@math.lu.se}

\address
{Department of Mathematics \& Computer Science, University of
Southern Denmark, Campusvej 55, DK-5230 Odense M, Denmark}
\email{svensson@imada.sdu.dk}

\begin{abstract}
We consider 4-dimensional Lie groups with left-invariant
Riemannian metrics.  For such groups we classify left-invariant
conformal foliations with minimal leaves of codimension 2.  These
foliations produce local complex-valued harmonic morphisms.
\end{abstract}

\maketitle

\section{Introduction}

This article develops an interplay between homogeneous spaces and
extremal mappings rather along the lines of Kaluza-Klein theory,
whereby the gravitational field (the metric) and electromagnetism
are unified in terms of a projection from a higher dimensional
space to a lower dimensional one. Although here, the dimensions
are not as in classical KK theory, see \cite{Bai-Lov}.

We are interested in the existence of complex-valued harmonic morphisms
from 4-dimensional Riemannian homogeneous spaces.  S. Ishihara has
in \cite{S-Ish} shown that any such space is either symmetric or
a Lie group equipped with a left-invariant metric, see also \cite{Ber}.
It is well-known that any 4-dimensional symmetric
space carries local harmonic morphisms and even global solutions
exist if the space is of non-compact type, see \cite{Gud-Sve-4}.
This means that we can focus our attention on the Lie group case.

The current authors have in \cite{Gud-Sve-5} {\it classified} the
3-dimensional Riemannian Lie groups $G$ carrying complex-valued
harmonic morphisms.  The classification is based on the fact that any
such local solution induces a global left-invariant conformal foliation
on $G$ with minimal leaves of codimension 2.

In \cite{Gud-Nor-1}, the authors introduce a general method for producing
left-invariant conformal foliations, on higher dimensional Riemannian
Lie groups, with minimal leaves.  This is then used to construct many
new examples, in particular, in the case when the the leaves are of
codimension 2.   This is important since such foliations produce local
complex-valued harmonic morphisms.

J. Nordstr\" om considers in \cite{Nor-1} a family of homogeneous
Hadamard manifolds of any dimension greater than 2.  They are all
Riemannian Lie groups with a rather simple bracket structure.  He
proves that none of these groups carry a left-invariant conformal
foliation with minimal fibres of codimension 2.  This makes it
interesting to understand what algebraic and geometric conditions
are necessary or sufficient for existence.

In the first part of this paper, we {\it classify} the left-invariant
conformal foliations on 4-dimensional Riemannian Lie groups, with
minimal leaves of codimension 2.  We prove the following result.

\begin{theorem}
Let $G$ be a 4-dimensional Lie group equipped with a left-invariant
Riemannian metric.  Let $\F$ be a left-invariant conformal foliation
on $G$ with minimal leaves of codimension 2.  Then the Lie algebra $\g$ of
$G$ belongs to one of the families $\g_1,\dots ,\g_{20}$ given below
and $\F$ is the corresponding foliation generated by $\g$.
\end{theorem}

We show that most of the complex-valued harmonic morphisms constructed
in this paper are {\it not} holomorphic with respect to any (integrable)
Hermitian structure on their domains.

In the second part, we consider Riemannian Lie groups of
higher dimensions.  On those we construct new conformal foliations
with minimal leaves of codimension 2.  Our examples show, for the
first time, that Theorem \ref{theo-Ric}
does not hold if the dimension of $M$ is greater than $3$.

\begin{theorem}\cite{Bai-Woo-1}\label{theo-Ric}
Let $(M,g)$ be a $3$-dimensional Riemannian manifold with a conformal
foliation $\F$ with minimal leaves of codimension $2$. Let $\{X,Y\}$
be a local orthonormal frame for the horizontal distribution $\H$.
Then the Ricci curvature satisfies
$$\text{Ric}(X,X)=\text{Ric}(Y,Y)\ \ \text{and}\ \ \text{Ric}(X,Y)=0.$$
\end{theorem}

This once more makes pertinent the connection with Kaluza-Klein
theories, where the Ricci tensor plays a crucial role.

For the general theory of harmonic morphisms between Riemannian
manifolds we refer to the excellent book \cite{Bai-Woo-book}
and the regularly updated on-line bibliography \cite{Gud-bib}.

\section{Harmonic morphisms and minimal conformal foliations}

Let $M$ and $N$ be two manifolds of dimensions $m$ and $n$,
respectively. A Riemannian metric $g$ on $M$ gives rise to the
notion of a {\it Laplacian} on $(M,g)$ and real-valued {\it
harmonic functions} $f:(M,g)\to\rn$. This can be generalized to
the concept of {\it harmonic maps} $\phi:(M,g)\to (N,h)$ between
Riemannian manifolds, which are solutions to a semi-linear system
of partial differential equations, see \cite{Bai-Woo-book}.

\begin{definition}
  A map $\phi:(M,g)\to (N,h)$ between Riemannian manifolds is
  called a {\it harmonic morphism} if, for any harmonic function
  $f:U\to\rn$ defined on an open subset $U$ of $N$ with $\phi^{-1}(U)$
non-empty,
  $f\circ\phi:\phi^{-1}(U)\to\rn$ is a harmonic function.
\end{definition}

The following characterization of harmonic morphisms between
Riemannian manifolds is due to Fuglede and T. Ishihara.  For the
definition of horizontal (weak) conformality we refer to
\cite{Bai-Woo-book}.

\begin{theorem}\cite{Fug-1,T-Ish}
  A map $\phi:(M,g)\to (N,h)$ between Riemannian manifolds is a
  harmonic morphism if and only if it is a horizontally (weakly)
  conformal harmonic map.
\end{theorem}

Let $(M,g)$ be a Riemannian manifold, $\V$ be an involutive
distribution on $M$ and denote by $\H$ its orthogonal
complement distribution on $M$.
As customary, we also use $\V$ and $\H$ to denote the
orthogonal projections onto the corresponding subbundles of $TM$
and denote by $\F$ the foliation tangent to
$\V$. The second fundamental form for $\V$ is given by
$$B^\V(U,V)=\frac 12\H(\nabla_UV+\nabla_VU)\qquad(U,V\in\V),$$
while the second fundamental form for $\H$ is given by
$$B^\H(X,Y)=\frac{1}{2}\V(\nabla_XY+\nabla_YX)\qquad(X,Y\in\H).$$
The foliation $\F$ tangent to $\V$ is said to be {\it conformal} if there is a
vector field $V\in \V$ such that $$B^\H=g\otimes V,$$ and
$\F$ is said to be {\it Riemannian} if $V=0$.
Furthermore, $\F$ is said to be {\it minimal} if $\text{trace}\ B^\V=0$ and
{\it totally geodesic} if $B^\V=0$. This is equivalent to the
leaves of $\F$ being minimal and totally geodesic submanifolds
of $M$, respectively.

It is easy to see that the fibres of a horizontally conformal
map (resp.\ Riemannian submersion) give rise to a conformal foliation
(resp.\ Riemannian foliation). Conversely, the leaves of any
conformal foliation (resp.\ Riemannian foliation) are
locally the fibres of a horizontally conformal map
(resp.\ Riemannian submersion), see \cite{Bai-Woo-book}.

The next result of Baird and Eells gives the theory of
harmonic morphisms, with values in a surface,
a strong geometric flavour.

\begin{theorem}\cite{Bai-Eel}\label{theo:B-E}
Let $\phi:(M^m,g)\to (N^2,h)$ be a horizontally conformal
submersion from a Riemannian manifold to a surface. Then $\phi$ is
harmonic if and only if $\phi$ has minimal fibres.
\end{theorem}

\section{4-dimensional Lie groups}

Let $G$ be a 4-dimensional Lie group equipped with a left-invariant
Riemannian metric.  Let $\g$ be the Lie algebra of $G$ and
$\{X,Y,Z,W\}$ be an orthonormal basis for $\g$.  Let $Z,W\in\g$
generate a 2-dimensional left-invariant and integrable distribution $\V$
on $G$ which is conformal and with minimal leaves.  We denote
by $\H$ the horizontal distribution, orthogonal to $\V$, generated by
$X,Y\in\g$.  Then it is easily seen that the
Lie bracket relations for $\g$ are of the form
\begin{eqnarray*}
\lb WZ&=&\lambda W,\\
\lb ZX&=&\alpha X +\beta Y+z_1 Z+w_1 W,\\
\lb ZY&=&-\beta X+\alpha Y+z_2 Z+w_2 W,\\
\lb WX&=&     a X     +b Y+z_3 Z-z_1W,\\
\lb WY&=&    -b X     +a Y+z_4 Z-z_2W,\\
\lb YX&=&     r X         +\theta_1 Z+\theta_2 W
\end{eqnarray*}
with real structure constants.  For later reference we state the following
easy result describing the geometry of the situation.

\begin{proposition}\label{prop-geometry}
Let $G$ be a 4-dimensional Lie group and $\{X,Y,Z,W\}$ be an
orthonormal basis for its Lie algebra as  above.  Then
\begin{enumerate}
\item[(i)] $\F$ is {\it totally geodesic} if and only if
  $z_1=z_2=z_3+w_1=z_4+w_2=0$,
\item[(ii)] $\F$ is {\it Riemannian} if and only if $\alpha=a=0$, and
\item[(iii)] $\H$ is {\it integrable} if and only if $\theta_1=\theta_2=0$.
\end{enumerate}
\end{proposition}

On the Riemannian Lie group $(G,g)$ there exist, up to sign,
exactly two invariant almost Hermitian structure $J_1$ and $J_2$
which are adapted to the orthogonal decomposition $\g=\V\oplus\H$
of the Lie algebra $\g$.  They are determined by
$$J_1X=Y,\ J_1Y=-X,\ J_1Z=W,\ J_1W=-Z,$$
$$J_2X=Y,\ J_2Y=-X,\ J_2W=Z,\ J_2Z=-W.$$
An elementary calculation involving the Nijenhuis tensor shows that
$J_1$ is integrable if and only if
$$2z_1-z_4-w_2=2z_2+z_3+w_1=0$$ and the same applies to $J_2$
if and only if
$$2z_1+z_4+w_2=2z_2-z_3-w_1=0.$$
This means that most of the complex-valued harmonic morphisms
constructed in this paper are {\it not} holomorphic with respect to any
Hermitian structure on the corresponding Lie groups.

\begin{remark}
This is interesting in the light of a result of J. C. Wood,
see \cite{Woo}. He shows that a submersive harmonic morphism from an
orientable 4-dimensional Einstein manifold $M^4$ to a Riemann surface,
or a conformal foliation of $M^4$ by minimal surfaces, determines an
(integrable) Hermitian structure with respect to which it is holomorphic.
\end{remark}

For any $V\in\V$, the adjoint action of
$V$ on $\H$ is conformal i.e.
$$\ip{\ad_VX}{Y}+\ip{X}{\ad_VY}=\rho\cdot\ip{X}{Y}\qquad(X,Y\in\H)$$
or put differently
$$\H\ad_V\big\vert_\H\in\rn\cdot\mathrm{Id}_\H+\so{\H}=\co{\H}.$$
Note that this is indeed a Lie algebra representation of $\V$, so that
$$\H\ad_{[Z,W]}\big\vert_\H
=[\H\ad_Z\big\vert_\H,\H\ad_W\big\vert_\H]\qquad(Z,W\in\V),$$
where the bracket on the right-hand side is just the usual
bracket on the space of endomorphisms on $\H$. This follows from the Jacobi
identity, and from the fact that $\V$ is integrable. Now, since
$\co{\H}$ is abelian, we see from this formula, that the adjoint
action of $[\V,\V]$ has no $\H$-component. This means that our analysis
branches into two different cases depending on whether $\V$ is abelian
or not.

\section{The case of non-abelian vertical distribution ($\lambda\neq 0$)}

We now assume that the vertical distribution $\V$ generated
by the left-invariant vector fields $Z,W\in\g$ is not abelian.
Then there exists an
orthonormal basis $\{Z,W\}$ for $\V$ such that $$[W,Z]=\lambda W$$
for some real constant $\lambda\neq 0$ and $\ad_W$ has no $\H$-component.
These conditions give the following bracket relations
\begin{eqnarray*}
\lb WZ&=&\lambda W,\\
\lb ZX&=&\alpha X+\beta Y+z_1Z+w_1 W,\\
\lb ZY&=&-\beta X+\alpha Y+z_2Z+w_2 W,\\
\lb WX&=&z_3Z-z_1W,\\
\lb WY&=&z_4Z-z_2W,\\
\lb YX&=&r X+\theta_1Z+\theta_2W.
\end{eqnarray*}
An elementary calculation shows that the two Jacobi equations
\begin{equation}\label{equa:Jacobi-1}
[[W,Z],X]+[[X,W],Z]+[[Z,X],W]=0,
\end{equation}
\begin{equation}\label{equa:Jacobi-2}
[[W,Z],Y]+[[Y,W],Z]+[[Z,Y],W]=0
\end{equation}
are equivalent to the following relations for the real
structure constants
\begin{equation}\label{equa-AB}
\begin{pmatrix}
\beta & \lambda-\alpha \\
\lambda-\alpha & -\beta
\end{pmatrix}
\begin{pmatrix}
z_1 & z_4 \\
z_2 & -z_3
\end{pmatrix}=
\begin{pmatrix}
0 & 0 \\
0 & 0
\end{pmatrix}.
\end{equation}

\subsection{Case (A) - ($\lambda\neq 0$ and $(\lambda-\alpha)^2+\beta^2\neq 0$)}

Applying equation (\ref{equa-AB}) we see that $z=0$ so the Lie brackets
satisfy
\begin{eqnarray*}
\lb WZ&=&\lambda W,\\
\lb ZX&=&\alpha X+\beta Y+w_1 W,\\
\lb ZY&=&-\beta X+\alpha Y+w_2 W,\\
\lb YX&=&r X+\theta_1Z+\theta_2W.
\end{eqnarray*}
In this situation it is easily seen that the Jacobi equations
\begin{equation}\label{equa:Jacobi-3}
[[X,Y],Z]+[[Z,X],Y]+[[Y,Z],X]=0,
\end{equation}
\begin{equation}\label{equa:Jacobi-4}
[[X,Y],W]+[[W,X],Y]+[[Y,W],X]=0
\end{equation}
are equivalent to
\begin{equation}
\theta_1=r\alpha=r\beta=0 \ \ \text{and}\ \ \theta_2(\lambda+2\alpha)=rw_1.
\end{equation}

\begin{example}[$\g_{1}(\lambda,r,w_1,w_2)$]\label{exam-A1}
If $r\neq 0$ then clearly $\alpha=\beta=0$ and $rw_1=\lambda\theta_2$.
This gives a 4-dimensional family of solutions satisfying
the following Lie bracket relations
\begin{eqnarray*}
\lb WZ&=&\lambda W,\\
\lb ZX&=&w_1 W,\\
\lb ZY&=&w_2 W,\\
\lambda \lb YX&=&\lambda r X+rw_1W.
\end{eqnarray*}
\end{example}

On the other hand, if $r=0$ then clearly
$\theta_1=\theta_2(\lambda+2\alpha)=0$
providing us with the following two examples.

\begin{example}[$\g_{2}(\lambda,\alpha,\beta,w_1,w_2)$]\label{exam-A2}
For $r=\theta_1=\theta_2=0$ we have the family
$\g=\g_2(\lambda,\alpha,\beta,w_1,w_2)$ given by
\begin{eqnarray*}
\lb WZ&=&\lambda W,\\
\lb ZX&=&\alpha X+\beta Y+w_1 W,\\
\lb ZY&=&-\beta X+\alpha Y+w_2 W.
\end{eqnarray*}
\end{example}

\begin{example}[$\g_{3}(\alpha,\beta,w_1,w_2,\theta_2)$]\label{exam-A3}
If $r=\theta_1=0$ and $\theta_2\neq 0$ then $\lambda=-2\alpha$ provides
us with the family
$\g=\g_3(\alpha,\beta,w_1,w_2,\theta_2)$ of solutions satisfying
\begin{eqnarray*}
\lb WZ&=&-2\alpha W,\\
\lb ZX&=&\alpha X+\beta  Y+w_1 W,\\
\lb ZY&=&-\beta  X+\alpha Y+w_2 W,\\
\lb YX&=&\theta_2W.
\end{eqnarray*}
\end{example}

\subsection{Case (B) - ($\lambda\neq 0$ and $(\lambda-\alpha)^2+\beta^2=0$)}

Under the assumptions that $\alpha=\lambda$ and $\beta=0$
we have the following bracket relations
\begin{eqnarray*}
\lb WZ&=&\lambda W,\\
\lb ZX&=&\lambda X+z_1Z+w_1 W,\\
\lb ZY&=&\lambda Y+z_2Z+w_2 W,\\
\lb WX&=&z_3Z-z_1W,\\
\lb WY&=&z_4Z-z_2W,\\
\lb YX&=&r X+\theta_1Z+\theta_2W.
\end{eqnarray*}
The Jacobi equations (\ref{equa:Jacobi-3}) and
(\ref{equa:Jacobi-4}) are easily seen to be equivalent to
$$z_1=z_3=z_4=\theta_1=0,\ \ z_2=-r\ \
\text{and}\ \ \lambda\theta_2=-z_2w_1.$$

\begin{example}[$\g_{4}(\lambda,z_2,w_1,w_2)$]\label{exam-B1}
In this case we have the family of solutions given by
\begin{eqnarray*}
\lb WZ&=&\lambda W,\\
\lb ZX&=&\lambda X+w_1 W,\\
\lb ZY&=&\lambda Y+z_2Z+w_2 W,\\
\lb WY&=&-z_2W,\\
\lambda\lb YX&=&-z_2\lambda X-z_2w_1W.
\end{eqnarray*}
\end{example}

\section{The case of abelian vertical distribution ($\lambda = 0$)}

We now assume that the vertical distribution $\V$ generated
by the left-invariant vector fields $Z,W\in\g$ is abelian.  Under this
natural algebraic condition, we have the general Lie bracket relations
\begin{eqnarray*}
\lb ZX&=&\alpha X +\beta Y+z_1 Z+w_1 W,\\
\lb ZY&=&-\beta X+\alpha Y+z_2 Z+w_2 W,\\
\lb WX&=&     a X     +b Y+z_3 Z-z_1W,\\
\lb WY&=&    -b X     +a Y+z_4 Z-z_2W,\\
\lb YX&=&     r X     +\theta_1 Z+\theta_2 W.
\end{eqnarray*}
An elementary calculation shows that the Jacobi equations
(\ref{equa:Jacobi-1}), (\ref{equa:Jacobi-2}),
(\ref{equa:Jacobi-3}) and (\ref{equa:Jacobi-4})
are equivalent to the following homogeneous systems of quadratic relations
for the real structure constants
\begin{equation}\label{equa-top}
-2\begin{pmatrix}
z_1 & w_1 \\
z_2 & w_2 \\
z_3 &-z_1 \\
z_4 &-z_2
\end{pmatrix}
\begin{pmatrix}
\alpha & \beta \\
a & b
\end{pmatrix}
=r
\begin{pmatrix}
 \beta & \alpha \\
\alpha & -\beta \\
b & a \\
a & -b
\end{pmatrix}
\end{equation}
and
\begin{equation}\label{equa-funny}
\begin{pmatrix}
-a\theta_2  \\
z_3w_2-z_4w_1 \\
2z_2z_3-2z_1z_4 \\
2z_1w_2-2z_2w_1
\end{pmatrix}
=
\begin{pmatrix}
\alpha\theta_1 \\
2a\theta_2+rz_1 \\
2a\theta_1-rz_3 \\
2\alpha\theta_2-rw_1
\end{pmatrix}.
\end{equation}

\subsection{Case (C) - ($\lambda=0$, $r\neq 0$ and $(a\beta-\alpha b)\neq 0$)}

According to equation (\ref{equa-top}) we can now solve $(z,w)$ in terms
of $(\alpha,\beta,a,b)$ and get
\begin{eqnarray}\label{equa-zw}
\begin{pmatrix}
z_1 & w_1 \\
z_2 & w_2 \\
z_3 &-z_1 \\
z_4 &-z_2
\end{pmatrix}
&=&\frac r{2(a\beta-\alpha b)}
\begin{pmatrix}
\beta b-\alpha a & \alpha^2-\beta^2 \\
\alpha b+\beta a & -2\alpha\beta \\
b^2-a^2 & \alpha a-\beta b \\
2ab & -(\alpha b+\beta a)
\end{pmatrix}.
\end{eqnarray}
By substituting these expressions for $(z,w)$ into equation
(\ref{equa-funny}) we then obtain the useful relations
\begin{equation}\label{equation-theta}
2(a\beta-\alpha b)
\begin{pmatrix}
a\theta_2 \\
a\theta_1 \\
\alpha\theta_2
\end{pmatrix}=r^2
\begin{pmatrix}
a\alpha\\
-a^2\\
\alpha^2
\end{pmatrix}\ \ \text{and}\ \ \alpha\theta_1=-a\theta_2.
\end{equation}
By solving those we get
\begin{equation*}
\theta_1=\frac{-ar^2}{2(a\beta-\alpha b)}\ \ \text{and}\ \
\theta_2= \frac{\alpha r^2}{2(a\beta-\alpha b)}.
\end{equation*}

\begin{example}[$\g_{5}(\alpha,a,\beta,b,r)$]\label{exam-C1}
In this case the Lie bracket relations take the following form
\begin{eqnarray*}
\lb ZX&=&\alpha X +\beta Y
+\frac{r(\beta b-\alpha a)}{2(a\beta-\alpha b)} Z+\frac{r(\alpha^2-\beta^2)}{2(a\beta-\alpha b)} W,\\
\lb ZY&=&-\beta X+\alpha Y
+\frac{r(\alpha b+\beta a)}{2(a\beta-\alpha b)} Z-\frac{r\alpha\beta}{(a\beta-\alpha b)} W,\\
\lb WX&=&     a X     +b Y
+\frac{r(b^2-a^2)}{2(a\beta-\alpha b)} Z+\frac{r(\alpha a-\beta b)}{2(a\beta-\alpha b)}W,\\
\lb WY&=&    -b X     +a Y
+\frac{rab }{(a\beta-\alpha b)} Z-\frac{r(\alpha b+\beta a)}{2(a\beta-\alpha b)}W,\\
\lb YX&=&     r X
-\frac{ar^2}{2(a\beta-\alpha b)} Z+\frac{\alpha r^2}{2(a\beta-\alpha b)} W.
\end{eqnarray*}
\end{example}

\subsection{Case (D) - ($\lambda=0$, $r\neq 0$ and $(a\beta-\alpha b)=0$)}

Picking the appropriate determinants from equation (\ref{equa-top})
we see that
\begin{equation*}
r^2(\alpha^2+\beta^2)=4(z_2w_1-z_1w_2)(\alpha b-a\beta)=0,
\end{equation*}
\begin{equation*}
r^2(a^2+b^2)=4(z_2z_3-z_1z_4)(\alpha b-a\beta)=0.
\end{equation*}
This means that in this case we have $\alpha=a=\beta=b=0$
so the corresponding foliations are all Riemannian.
The equation (\ref{equa-top}) is  automatically satisfied
and the system (\ref{equa-funny}) takes the form
\begin{equation}\label{equa-funny-simple}
\begin{pmatrix}
z_3w_2-z_4w_1 \\
2z_2z_3-2z_1z_4 \\
2z_1w_2-2z_2w_1
\end{pmatrix}
=r
\begin{pmatrix}
z_1 \\
-z_3 \\
-w_1
\end{pmatrix}.
\end{equation}
Applying basic algebraic manipulations on (\ref{equa-funny-simple}) we obtain
\begin{equation}\label{equa-D}
z_1^2+w_1z_3=0\ \ \text{and}\ \ 2z_1z_2+z_4w_1+z_3w_2=0.
\end{equation}

\begin{example}[$\g_{6}(z_1,z_2,z_3,r,\theta_1,\theta_2)$]\label{exam-D1}
When $z_1^2=-w_1z_3\neq 0$ it immediately follows that
$$z_4 =\frac{z_3(r+2z_2)}{2z_1},\ \ w_1 = -\frac{z_1^2}{z_3},
\ \ w_2 = \frac{z_1(r-2z_2)}{2z_3}.$$
This provides the following family of solutions
\begin{eqnarray*}
\lb ZX&=&z_1 Z-\frac{z_1^2}{z_3} W,\\
\lb ZY&=&z_2 Z+\frac{z_1(r-2z_2)}{2z_3} W,\\
\lb WX&=&z_3 Z-z_1W,\\
\lb WY&=&\frac{z_3(r+2z_2)}{2z_1} Z-z_2W,\\
\lb YX&=&r X+\theta_1 Z+\theta_2 W.
\end{eqnarray*}
\end{example}

When $z_1=0$ the system (\ref{equa-funny-simple}) is equivalent to
\begin{equation}\label{equa-funny-too}
z_3w_2=z_4w_1,\ \ z_3(2z_2+r)=0,\ \ w_1(2z_2-r)=0.
\end{equation}

\begin{example}[$\g_{7}(z_2,w_1,w_2,\theta_1,\theta_2)$]\label{exam-D2}
When $z_1=0$ and $w_1\neq 0$ it follows from (\ref{equa-funny-too})
that $z_3=z_4=0$ and $r=2z_2$ so we have the following solutions
\begin{eqnarray*}
\lb ZX&=&w_1 W,\\
\lb ZY&=&z_2 Z+w_2 W,\\
\lb WY&=&-z_2 W,\\
\lb YX&=&2z_2 X+\theta_1 Z+\theta_2 W.
\end{eqnarray*}
\end{example}

When $z_1=w_1=0$ it follows that $z_3w_2=0$ and $z_3(2z_2+r)=0$.
The two possible cases are given in Examples \ref{exam-D3} and \ref{exam-D4}.

\begin{example}[$\g_{8}(z_2,z_4,w_2,r,\theta_1,\theta_2)$]\label{exam-D3}
When $z_1=z_3=w_1=0$ we have the solutions
\begin{eqnarray*}
\lb ZY&=&z_2 Z+w_2 W,\\
\lb WY&=&z_4 Z-z_2 W,\\
\lb YX&=&r X+\theta_1 Z+\theta_2 W.
\end{eqnarray*}
\end{example}

\begin{example}[$\g_{9}(z_2,z_3,z_4,\theta_1,\theta_2)$]\label{exam-D4}
The conditions $z_1=w_1=w_2=0$ and $r=-2z_2$ produce the following solutions
\begin{eqnarray*}
\lb ZY&=&z_2 Z,\\
\lb WX&=&z_3 Z,\\
\lb WY&=&z_4 Z-z_2 W,\\
\lb YX&=&-2z_2 X+\theta_1 Z+\theta_2 W.
\end{eqnarray*}
\end{example}

\subsection{Case (E) - ($\lambda=0$, $r=0$ and $\alpha b-a\beta\neq 0$)}

In this situation the equation (\ref{equa-top}) takes the following simple form
\begin{equation*}
\begin{pmatrix}
z_1 & w_1 \\
z_2 & w_2 \\
z_3 &-z_1 \\
z_4 &-z_2
\end{pmatrix}
\begin{pmatrix}
\alpha & \beta \\
a & b
\end{pmatrix}
=
\begin{pmatrix}
0 & 0 \\
0 & 0 \\
0 & 0 \\
0 & 0
\end{pmatrix}.
\end{equation*}
As an immediately consequence of $\alpha b-a\beta\neq 0$ we see that $z=w=0$.
This means that the system (\ref{equa-funny}) is equivalent to
$$\alpha\theta_1=a\theta_2=a\theta_1=\alpha\theta_2=0.$$
Applying $\alpha b\neq a\beta$ we see that $\theta_1=\theta_2=0$ so we
have the following.

\begin{example}[$\g_{10}(\alpha,a,\beta,b)$]\label{exam-E1}
If $\alpha b-a\beta\neq 0$ then the solutions are given by
\begin{eqnarray*}
\lb ZX&=&\alpha X+ \beta Y,\\
\lb ZY&=&-\beta X+\alpha Y,\\
\lb WX&=&     a X     +b Y,\\
\lb WY&=&    -b X     +a Y.
\end{eqnarray*}
\end{example}

\subsection{Case (F) - ($\lambda=0$, $r=0$ and $\alpha b-a\beta= 0$)}

In this section we assume that both $\lambda$ and $r$ are zero. This
means that the special choice of bases $\{X,Y\}$ and $\{Z,Y\}$ for
$\H$ and $\V$, we made above, is irrelevant and that we get the
following symmetric system of bracket relations
\begin{eqnarray}\label{symmetric-system}
\lb ZX&=&\alpha X +\beta Y+z_1 Z+w_1 W,\\
\lb ZY&=&-\beta X+\alpha Y+z_2 Z+w_2 W,\\
\lb WX&=&     a X     +b Y+z_3 Z-z_1W,\\
\lb WY&=&    -b X     +a Y+z_4 Z-z_2W,\\
\lb YX&=&\theta_1 Z+\theta_2 W.
\end{eqnarray}
The Jacobi identities now take the following form
\begin{equation}\label{Jacobi-D2}
\begin{pmatrix}
\alpha z_1 & \beta z_1 \\
\alpha z_2 & \beta z_2 \\
\alpha z_3 & \beta z_3 \\
\alpha z_4 & \beta z_4
\end{pmatrix}
=
\begin{pmatrix}
-aw_1 & -bw_1 \\
-aw_2 & -bw_2 \\
az_1 &bz_1 \\
az_2 &bz_2
\end{pmatrix}
,\ \
\begin{pmatrix}
-a\theta_2  \\
z_3w_2-z_4w_1 \\
z_2z_3-z_1z_4 \\
z_1w_2-z_2w_1
\end{pmatrix}=
\begin{pmatrix}
\alpha\theta_1 \\
2a\theta_2 \\
a\theta_1 \\
\alpha\theta_2
\end{pmatrix}.
\end{equation}

Our analysis divides into disjoint cases parametrized by
$\Lambda=(\alpha,a,\beta,b)$.  The variables are assumed to be zero
if and only if they are marked by $0$.
For example, if $\Lambda=(0,a,\beta,0)$ then the two variable $\alpha$
and $b$ are assumed to be zero and $a$ and $\beta$ to be non-zero.

Let us first consider the case when $\Lambda=(0,0,0,0)$.
Then the Jacobi equations are equivalent to
\begin{equation*}
\begin{pmatrix}
z_3w_2-z_4w_1 \\
z_2z_3-z_1z_4 \\
z_1w_2-z_2w_1
\end{pmatrix}=
\begin{pmatrix}
0 \\
0 \\
0
\end{pmatrix}
\end{equation*}
and all the possible cases are covered by Examples \ref{exam-F1}-\ref{exam-F4}.

\begin{example}[$\g_{11}(z_1,z_2,z_3,w_1,\theta_1,\theta_2)$]\label{exam-F1}
If $\Lambda=(0,0,0,0)$ and $z_1\neq 0$
then we have the following solutions
\begin{eqnarray*}
\lb ZX&=&z_1 Z+w_1 W,\\
\lb ZY&=&z_2 Z+\frac{z_2w_1}{z_1} W,\\
\lb WX&=&z_3 Z-z_1W,\\
\lb WY&=&\frac{z_2z_3}{z_1} Z-z_2 W,\\
\lb YX&=&\theta_1 Z+\theta_2 W.
\end{eqnarray*}
\end{example}

If $z_1=0$ then $$z_3w_2=z_4w_1,\ \ z_2z_3=0,\ \ z_2w_1=0.$$

\begin{example}[$\g_{12}(z_3,w_1,w_2,\theta_1,\theta_2)$]\label{exam-F2}
If $\Lambda=(0,0,0,0)$, $z_1=0$ and $w_1\neq 0$ then
$$z_2=0\ \ \text{and}\ \ z_4=\frac{z_3w_2}{w_1}$$
so possible solutions are given by
\begin{eqnarray*}
\lb ZX&=&w_1 W,\\
\lb ZY&=&w_2 W,\\
\lb WX&=&z_3 Z,\\
\lb WY&=&\frac{z_3w_2}{w_1} Z,\\
\lb YX&=&\theta_1 Z+\theta_2 W.
\end{eqnarray*}
\end{example}

\begin{example}[$\g_{13}(z_3,z_4,\theta_1,\theta_2)$]\label{exam-F3}
If $\Lambda=(0,0,0,0)$,  $z_1=w_1=0$ and $z_3\neq 0$ then $z_2=w_2=0$
producing the solutions
\begin{eqnarray*}
\lb WX&=&z_3 Z,\\
\lb WY&=&z_4 Z,\\
\lb YX&=&\theta_1 Z+\theta_2 W.
\end{eqnarray*}
\end{example}

\begin{example}[$\g_{14}(z_2,z_4,w_2,\theta_1,\theta_2)$]\label{exam-F4}
The conditions $\Lambda=(0,0,0,0)$ and $z_1=z_3=w_1=0$
provide the solutions
\begin{eqnarray*}
\lb ZY&=&z_2 Z+w_2 W,\\
\lb WY&=&z_4 Z-z_2 W,\\
\lb YX&=&\theta_1 Z+\theta_2 W.
\end{eqnarray*}
\end{example}

The rest of this section is devoted to the situation
$\Lambda\neq(0,0,0,0)$.  Up to obvious similarities,
all the possible cases are listed in Examples
\ref{exam-F5}-\ref{exam-F10}.

\begin{example}[$\g_{15}(\alpha,w_1,w_2)$]\label{exam-F5}
For $\Lambda=(\alpha,0,0,0)$ the Jacobi identities (\ref{Jacobi-D2}) give
\begin{equation*}
z_1=z_2=z_3=z_4=\theta_1=\theta_2=0.
\end{equation*}
This yields the 3-dimensional family of solutions
\begin{eqnarray*}
\lb ZX&=&\alpha X+w_1 W,\\
\lb ZY&=&\alpha Y+w_2 W.
\end{eqnarray*}
It follows from the symmetry of (\ref{symmetric-system})
that $\Lambda=(0,a,0,0)$ gives similar solutions.
\end{example}

\begin{example}[$\g_{16}(\beta,w_1,w_2,\theta_1,\theta_2)$]\label{exam-F6}
If $\Lambda=(0,0,\beta,0)$ then the Jacobi identities are equivalent to
\begin{equation*}
z_1=z_2=z_3=z_4=0.
\end{equation*}
This provides us with the 5-dimensional family of solutions satisfying
\begin{eqnarray*}
\lb ZX&=&\beta Y+w_1 W,\\
\lb ZY&=&-\beta X+w_2 W,\\
\lb YX&=&\theta_1 Z+\theta_2 W.
\end{eqnarray*}
The case $\Lambda=(0,0,0,b)$ is similar.
\end{example}

\begin{example}[$\g_{17}(\alpha,a,w_1,w_2)$]\label{exam-F7}
In the case when $\Lambda=(\alpha,a,0,0)$ the Jacobi equations give
$$
z_1=-\frac {aw_1}\alpha,
\ \ z_2=-\frac {aw_2}\alpha,
\ \ z_3=-\frac {a^2w_1}{\alpha^2},
\ \ z_4=-\frac {a^2w_2}{\alpha^2},\ \
\theta_1=0,\ \ \theta_2=0
$$
so we obtain the 4-dimensional family of solutions
\begin{eqnarray*}
\lb ZX&=&\alpha X-\frac {aw_1}\alpha Z+w_1 W,\\
\lb ZY&=&\alpha Y-\frac {aw_2}\alpha Z+w_2 W,\\
\lb WX&=&a X-\frac {a^2w_1}{\alpha^2} Z+\frac {aw_1}\alpha W,\\
\lb WY&=&a Y-\frac {a^2w_2}{\alpha^2} Z+\frac {aw_2}\alpha W.
\end{eqnarray*}
\end{example}

\begin{example}[$\g_{18}(\beta,b,z_3,z_4,\theta_1,\theta_2)$]\label{exam-F8}
In the case when $\Lambda=(0,0,\beta,b)$ the Jacobi equations give
$$z_1=\frac {\beta z_3}b,
\ \ z_2=\frac {\beta z_4}b,
\ \ w_1=-\frac {\beta^2z_3}{b^2},
\ \ w_2=-\frac {\beta^2z_4}{b^2}
$$
so we have the following family of solutions
\begin{eqnarray*}
\lb ZX&=& \beta Y+\frac {\beta z_3}b Z-\frac {\beta^2z_3}{b^2} W,\\
\lb ZY&=&-\beta X+\frac {\beta z_4}b Z-\frac {\beta^2z_4}{b^2} W,\\
\lb WX&=& b Y+z_3 Z-\frac {\beta z_3}b W,\\
\lb WY&=&-b X+z_4 Z-\frac {\beta z_4}b W,\\
\lb YX&=&\theta_1 Z+\theta_2 W.
\end{eqnarray*}
\end{example}

\begin{example}[$\g_{19}(\alpha,\beta,w_1,w_2)$]\label{exam-F9}
When $\Lambda=(\alpha,0,\beta,0)$ we get from the Jacobi equations that
$$z_1=z_2=z_3=z_4=\theta_1=\theta_2=0.$$
This produces the solutions
\begin{eqnarray*}
\lb ZX&=&\alpha X +\beta Y+w_1 W,\\
\lb ZY&=&-\beta X+\alpha Y+w_2 W.
\end{eqnarray*}
The case $\Lambda=(0,a,0,b)$ is similar.
\end{example}

\begin{example}[$\g_{20}(\alpha,a,\beta,w_1,w_2)$]\label{exam-F10}
For $\Lambda=(\alpha,a,\beta,b)$ the Jacobi relation give
$$
z_1=-\frac {aw_1}\alpha,
\ \ z_2=-\frac {aw_2}\alpha,
\ \ z_3=-\frac {a^2w_1}{\alpha^2},
\ \ z_4=-\frac {a^2w_2}{\alpha^2},$$
$$b=\frac{\beta a}\alpha,\ \ \theta_1=0,\ \ \theta_2=0.$$
In this case we have solutions of the form
\begin{eqnarray*}
\lb ZX&=&\alpha X+\beta Y-\frac {aw_1}\alpha Z+w_1 W,\\
\lb ZY&=&-\beta X+\alpha Y-\frac {aw_2}\alpha Z+w_2 W,\\
\lb WX&=&a X+\frac{\beta a}\alpha Y-\frac {a^2w_1}{\alpha^2} Z+\frac a\alpha w_1 W,\\
\lb WY&=&-\frac{\beta a}\alpha X+a Y-\frac {a^2w_2}{\alpha^2} Z+\frac a\alpha w_2 W.
\end{eqnarray*}
\end{example}

\section{The Ricci Operator}

In this section we produce Riemannian Lie groups carrying conformal
foliations with minimal leaves. We start with an infinite series of
Riemannian Lie groups which are all nilpotent.

\begin{example}
For a positive integer $n\in\zn^+$ let $\text{Nil}^{n+2}$
be the simply connected
nilpotent Lie group with Lie algebra $\n_{n+2}=\rn\ltimes\rn^{n+1}$
generated by the left-invariant
vector fields $W,X_1,\dots ,X_{n+1}$ satisfying
$$[W,X_k]=X_{k+1},\ \ k=1,2,\dots ,n.$$  Equip $\text{Nil}^{n+2}$
with the left-invariant
Riemannian metric such that $$\{W,X_1,\dots ,X_{n+1}\}$$ is an
orthonormal basis
for the Lie algebra $\n_{n+2}$.
It is easily shown that for each $k=1,2,\dots ,n$ the orthogonal
decomposition
$$\H_k\oplus\V_k=<W,X_1,\dots ,X_k>\oplus <X_{k+1},\dots ,X_{n+1}>$$
gives an example of a minimal Riemannian foliation, with
non-integrable horizontal
distribution and leaves which are not totally geodesic.

A standard calculation of the Ricci operator
$\text{Ric}:\n_{n+2}\to\n_{n+2}$ of
$\text{Nil}^{n+2}$ shows that
$$\text{Ric}(X_2)=\dots =\text{Ric}(X_n)=0,$$
$$\text{Ric}(X_1)=-\frac 12 X_1,
\ \ \text{Ric}(X_{n+1})=\frac 12X_{n+1},\ \ \text{Ric}(W)=-\frac n2W.$$
This implies that for the orthogonal decomposition
$$\H_1\oplus\V_1=<W,X_1>\oplus <X_2,\dots ,X_{n+1}>$$
the Ricci curvature tensor satisfies
$$\text{Ric}(X_1,X_1)-\text{Ric}(W,W)=\frac{n-1}2.$$

\end{example}

We complete this section with an infinite series of solvable
Riemannian Lie groups with conformal minimal foliations.

\begin{example}
Let $Sol^{n+1}$ be the simply connected Lie group with Lie algebra
$$\s_{n+1}=\rn\ltimes\rn^n$$
generated by the left-invariant vector fields $W,X_1,\dots ,X_{n}$
satisfying
$$[W,X_k]=\alpha_kX_{k},\ \ k=1,2,\dots ,n, \ \ \alpha_k\in\rn.$$
Then for each $k=1,2,\dots ,n-2$
the orthogonal decomposition
$$\H_k\oplus\V_k=<W,X_1,\dots ,X_k>\oplus <X_{k+1},\dots ,X_{n}>$$
gives a solution if and only if $$\alpha_{k+1}+\dots +\alpha_n=0.$$

A standard calculation of the Ricci operator
$\text{Ric}:\s_{n+1}\to\s_{n+1}$ of $\text{Sol}^{n+1}$
shows that
$$\text{Ric}(X_k)=-\alpha_k(\alpha_1+\dots +\alpha_n)X_k,
\ \ \text{Ric}(W)=-|\alpha |^2W.$$

If $k=1$ and $n\ge 3$ we can easily find $\alpha_1,\dots ,\alpha_n$
such that $\alpha_{2}+\dots +\alpha_n=0$ but
$\alpha_{2}^2+\dots +\alpha_n^2\neq 0$.  Again this shows that
Theorem \ref{theo-Ric} does not hold in any dimension greater than $3$.
\end{example}

\section{Acknowledgements}
The authors are grateful to Jonas Nordstr\" om for useful
discussions on this work.  They would also like
to thank the referee for pointing out the interesting
connection with Kaluza-Klein theory. 
The second author was supported by the Danish Council for
Independent Research under the project 
{\it Symmetry Techniques in Differential Geometry}.

\end{document}